\newif\ifjournal\journalfalse   

\ifjournal
\documentclass{elsart}
\usepackage{latexsym,amsmath,amssymb,hyperref}
\else 
\documentclass[12pt,twoside]{article}
\usepackage{latexsym,amsmath,amssymb,hyperref}
\fi
\ifjournal
\def\beq{\begin{equation}}    \def\eeq{\end{equation}}
\def\beqn{\begin{displaymath}}\def\eeqn{\end{displaymath}}
\def\bqa{\begin{eqnarray}}    \def\eqa{\end{eqnarray}}
\def\bqan{\begin{eqnarray*}}  \def\eqan{\end{eqnarray*}}
\else
\def\,{\mskip 3mu} \def\>{\mskip 4mu plus 2mu minus 4mu} \def\;{\mskip 5mu plus 5mu} \def\!{\mskip-3mu}
\def\dispmuskip{\thinmuskip= 3mu plus 0mu minus 2mu \medmuskip=  4mu plus 2mu minus 2mu \thickmuskip=5mu plus 5mu minus 2mu}
\def\textmuskip{\thinmuskip= 0mu                    \medmuskip=  1mu plus 1mu minus 1mu \thickmuskip=2mu plus 3mu minus 1mu}
\textmuskip
\def\beq{\dispmuskip\begin{equation}}    \def\eeq{\end{equation}\textmuskip}
\def\beqn{\dispmuskip\begin{displaymath}}\def\eeqn{\end{displaymath}\textmuskip}
\def\bqa{\dispmuskip\begin{eqnarray}}    \def\eqa{\end{eqnarray}\textmuskip}
\def\bqan{\dispmuskip\begin{eqnarray*}}  \def\eqan{\end{eqnarray*}\textmuskip}
\fi

\ifjournal\else
\newenvironment{keywords}{\centerline{\bf\small
Keywords}\begin{quote}\small}{\par\end{quote}\vskip 1ex}
\fi
\ifjournal
\def\paradot#1{{\itshape{#1.}}}

\else
\def\paradot#1{\vspace{1.3ex plus 0.5ex minus 0.5ex}\noindent{\bf\boldmath{#1.}}}

\fi
\def\aidx#1{}

\def\indxs#1#2{}
\def\req#1{(\ref{#1})}
\def\toinfty#1{\smash{\stackrel{#1\to\infty}{\longrightarrow}}}

\def\eps{\varepsilon}
\def\epstr{\epsilon}                    
\def\nq{\hspace{-1em}}
\def\qed{\hspace*{\fill}\rule{1.4ex}{1.4ex}$\quad$\\}
\def\odt{{\textstyle{1\over 2}}}
\def\odn{{\textstyle{1\over n}}}
\def\SetR{I\!\!R}\def\SetN{I\!\!N}\def\SetQ{I\!\!\!Q}
\def\M{{\cal M}}                        
\def\S{{\cal S}}                        
\def\X{{\cal X}}                        
\def\qmbox#1{{\quad\mbox{#1}\quad}}
\def\qqmbox#1{{\qquad\mbox{#1}\qquad}}
\def\equa{\,\smash{\stackrel{\raisebox{0.8ex}{$\scriptstyle+$}}{\smash=}}\,}
\def\leqa{\,\smash{\stackrel{\raisebox{1ex}{$\scriptstyle\!\!\;+$}}{\smash\leq}}\,}
\def\geqa{\,\smash{\stackrel{\raisebox{1ex}{$\scriptstyle+\!\!\;$}}{\smash\geq}}\,}
\def\eqm{\,\smash{\stackrel{\raisebox{0.6ex}{$\scriptstyle\times$}}{\smash=}}\,}
\def\leqm{\,\smash{\stackrel{\raisebox{1ex}{$\scriptstyle\!\!\;\times$}}{\smash\leq}}\,}
\def\geqm{\,\smash{\stackrel{\raisebox{1ex}{$\scriptstyle\times\!\!\;$}}{\smash\geq}}\,}
\def\th{\theta}
\def\e{{\rm e}}                        
\def\B{\{0,1\}}
\def\E{{\bf E}}                         
\def\P{{\rm P}}                         
\def\l{\ell}
\def\lb{{\log_2}}

\def\v{\boldsymbol}
\def\text#1{\mbox{\scriptsize #1}}

\def\o{\omega}
\def\a{\alpha}
\def\Km{K\!m}
\def\thb{{\th'}}
\def\thp{\th'}

\def\Thp{\Theta'}
\def\ir{\em}                   
\def\eqbr#1{\discretionary{#1}{#1}{#1}}

\begin{document}

\ifjournal

\begin{frontmatter}
\title{On Universal Prediction \\ and Bayesian Confirmation}
\author{Marcus Hutter}
\address{RSISE @ ANU and SML @ NICTA \\
Canberra, ACT, 0200, Australia \\
\texttt{marcus@hutter1.net \ \  www.hutter1.net} }

\else

\title{\vspace{-4ex}\normalsize
\vskip 2mm\bf\Large\hrule height5pt \vskip 2mm
On Universal Prediction \\ and Bayesian Confirmation
\vskip 2mm \hrule height2pt}
\author{{\bf Marcus Hutter}\\[3mm]
\normalsize RSISE$\,$@$\,$ANU and SML$\,$@$\,$NICTA \\
\normalsize Canberra, ACT, 0200, Australia \\
\normalsize \texttt{marcus@hutter1.net \ \  www.hutter1.net}}
\date{11 September 2007}
\maketitle
\vspace*{-4ex}

\fi

\begin{abstract}\noindent
The Bayesian framework is a well-studied and successful framework
for inductive reasoning, which includes hypothesis testing and
confirmation, parameter estimation, sequence prediction,
classification, and regression. But standard statistical guidelines
for choosing the model class and prior are not always available or
can fail, in particular in complex situations.
Solomonoff completed the Bayesian framework by providing a
rigorous, unique, formal, and universal choice for the
model class and the prior.
I discuss in breadth how and in which sense universal
(non-i.i.d.)\ sequence prediction solves various (philosophical)
problems of traditional Bayesian sequence prediction.
I show that Solomonoff's model possesses many desirable
properties: Strong total and future bounds, and weak
instantaneous bounds, and in contrast to most classical continuous
prior densities has no zero p(oste)rior problem, i.e.\ can confirm
universal hypotheses, is reparametrization and regrouping
invariant, and avoids the old-evidence and updating problem. It
even performs well (actually better) in non-computable
environments.
\ifjournal\else\vspace{1ex}
\def\contentsname{\centering\normalsize Contents}
{\parskip=-2.5ex\tableofcontents}\fi
\end{abstract}

\ifjournal\begin{keyword}\else\begin{keywords}\fi
Sequence prediction, %
Bayes, %
Solomonoff prior, %
Kolmogorov complexity, %
Occam's razor, %
prediction bounds, %
model classes, %
philosophical issues, %
symmetry principle, %
confirmation theory, %
Black raven paradox, %
reparametrization invariance, %
old-evidence/updating problem, %
(non)computable environments.
\ifjournal\end{keyword}\else\end{keywords}\fi

\ifjournal\end{frontmatter}\else\pagebreak\fi

\ifjournal\else\newpage\fi
\section{Introduction}\label{secInt}

\begin{quote}\it
``... in spite of it's incomputability, Algorithmic Probability can serve as a kind of
`Gold Standard' for induction systems'' \par
\hfill --- {\sl Ray Solomonoff (1997)}
\end{quote}

Given the weather in the past, what is the probability of rain
tomorrow? What is the correct answer in an IQ test asking to
continue the sequence 1,4,9,16,? Given historic stock-charts, can
one predict the quotes of tomorrow? Assuming the sun rose 5000
years every day, how likely is doomsday (that the sun does not
rise) tomorrow? These are instances of the important problem of
induction or time-series forecasting or sequence
prediction. Finding prediction rules for every particular (new)
problem is possible but cumbersome and prone to disagreement or
contradiction. What is desirable is a formal general
theory for prediction.

The Bayesian framework is the most consistent and successful
framework developed thus far \cite{Earman:93,Jaynes:03}. A Bayesian
considers a set of environments\eqbr=hypotheses\eqbr=models $\M$
which includes the true data generating probability distribution
$\mu$. From one's prior belief $w_\nu$ in environment $\nu\in\M$ and
the observed data sequence $x=x_1...x_n$, Bayes' rule yields one's
posterior confidence in $\nu$. In a prequential \cite{Dawid:84} or
transductive \cite[Sec.9.1]{Vapnik:99} setting, one directly
determines the predictive probability of the next symbol $x_{n+1}$
without the intermediate step of identifying a (true or good or
causal or useful) model. With the exception of Section \ref{secIID},
this paper concentrates on {\em prediction} rather than model
identification. The ultimate goal is to make ``good'' predictions in
the sense of maximizing one's profit or minimizing one's loss. Note
that classification and regression can be regarded as special
sequence prediction problems, where the sequence $x_1 y_1...x_n y_n
x_{n+1}$ of $(x,y)$-pairs is given and the class label or function
value $y_{n+1}$ shall be predicted.

The Bayesian framework leaves open how to choose the model class
$\M$ and prior $w_\nu$. General guidelines are that $\M$
should be small but large enough to contain the true environment $\mu$,
and $w_\nu$ should reflect one's prior (subjective) belief in $\nu$
or should be non-informative or neutral or objective if no prior
knowledge is available. But these are informal and ambiguous
considerations outside the formal Bayesian framework.
Solomonoff's \cite{Solomonoff:64} rigorous, essentially unique,
formal, and universal solution to this problem is to consider a
single large universal class $\M_U$ suitable for {\em all}
induction problems. The corresponding universal prior $w_\nu^U$ is
biased towards simple environments in such a way that it
dominates (=superior to) all other priors. This leads to an a
priori probability $M(x)$ which is equivalent to the probability
that a universal Turing machine with random input tape outputs
$x$, and the shortest program computing $x$ produces the most
likely continuation (prediction) of $x$.

Many interesting, important, and deep results have been proven for
Solomonoff's universal distribution $M$
\cite{Zvonkin:70,Solomonoff:78,Gacs:83,Li:97,Hutter:01errbnd,Hutter:04uaibook}.
The motivation and goal of this paper is %
to provide a broad discussion of how and in which sense universal sequence prediction
solves all kinds of (philosophical) problems
of Bayesian sequence prediction, and to %
present some recent results. %
Many arguments and ideas could be further developed. I hope that
the exposition stimulates such a future, more detailed, investigation.

In Section \ref{secBSP}, I review the excellent predictive and
decision-theoretic performance results of Bayesian sequence
prediction for generic (non-i.i.d.)\ countable and continuous model
classes.
Section \ref{secPrior} critically reviews the classical principles
(indifference, symmetry, minimax) for obtaining objective priors,
introduces the universal prior inspired by Occam's razor
and quantified in terms of Kolmogorov complexity.
In Section \ref{secIID} (for i.i.d.\ $\M$) and Section
\ref{secUSP} (for universal $\M_U$) I show various desirable
properties of the universal prior and class (non-zero
p(oste)rior, confirmation of universal hypotheses,
reparametrization and regrouping invariance, no old-evidence and
updating problem) in contrast to (most) classical continuous
prior densities. I also complement the general total bounds of
Section \ref{secBSP} with some universal and some i.i.d.-specific
instantaneous and future bounds. Finally, I
show that the universal mixture performs better than classical
continuous mixtures, even in uncomputable environments.
Section \ref{secDisc} contains critique, summary, and conclusions.

The reparametrization and regrouping invariance, the (weak)
instantaneous bounds, the good performance of $M$ in
non-computable environments, and most of the discussion (zero
prior and universal hypotheses, old evidence) are new or new in
the light of universal sequence prediction.
Technical and mathematical non-trivial new results are the
Hellinger-like loss bound \req{lbnd} and the instantaneous bounds
\req{iIIDbnd} and \req{iMbnd}.

\section{Bayesian Sequence Prediction}\label{secBSP}

I now formally introduce the Bayesian sequence prediction setup
and describe the most important results. I consider sequences
over a finite alphabet, assume that the true environment is
unknown but known to belong to a countable or continuous
class of environments (no i.i.d.\ or Markov or stationarity
assumption), and consider general prior. I show that the
predictive distribution converges rapidly to the true sampling
distribution and that the Bayes-optimal predictor performs
excellent for any bounded loss function.

\paradot{Notation}
I use letters $t,n\in\SetN$ for natural numbers, and
denote the cardinality of a set $\cal S$ by $\#{\cal S}$ or $|{\cal S}|$.
I write $\X^*$ for the set of finite strings over some alphabet
$\X$, and $\X^\infty$ for the set of infinite sequences.
For a string $x\in\X^*$ of length $\l(x)=n$ I write
$x_1x_2...x_n$ with $x_t\in\X$, and further abbreviate
$x_{t:n}:=x_t x_{t+1}...x_{n-1}x_n$ and $x_{<n}:=x_1... x_{n-1}$.

I assume that sequence $\o=\o_{1:\infty}\in\X^\infty$ is sampled
from the ``true'' probability measure $\mu$, i.e.\
$\mu(x_{1:n}):=\P[\o_{1:n}=x_{1:n}|\mu]$ is the $\mu$-probability
that $\o$ starts with $x_{1:n}$. I denote expectations w.r.t.\ $\mu$
by $\E$. In particular for a function $f:\X^n\to\SetR$, we have
$\E[f]=\E[f(\o_{1:n})]=\sum_{x_{1:n}}\mu(x_{1:n})f(x_{1:n})$.
Note that in Bayesian learning, measures, environments, and models
coincide, and are the same objects; let $\M=\{\nu_1,\nu_2,...\}$
denote a countable class of these measures. Assume that %
(a) $\mu$ is unknown but known to be a member of ${\cal M}$, %
(b) $\{H_\nu:\nu\in\M\}$ forms a mutually exclusive and complete
class of hypotheses, and (c) %
$w_\nu:=\P[H_\nu]$ is the given prior belief in $H_{\nu}$.
Then $\xi(x_{1:n}):=\P[\o_{1:n}=x_{1:n}] =
\sum_{\nu\in\M}\P[\o_{1:n}=x_{1:n}|H_\nu]\P[H_\nu]$ must be our
(prior) belief in $x_{1:n}$, and
$w_\nu(x_{1:n}):=\P[H_\nu|\o_{1:n}=x_{1:n}]=
{\P[\o_{1:n}=x_{1:n}|H_\nu]\P[H_\nu]\over\P[\o_{1:n}=x_{1:n}]}$ be
our posterior belief in $\nu$ by Bayes' rule.

For a sequence $a_1, a_2, ...$ of random variables,
$\sum_{t=1}^\infty\E[a_t^2]\leq c<\infty$ implies $a_t\toinfty{t} 0$
with $\mu$-probability 1 (w.p.1). Convergence is rapid in the
sense that the probability that $a_t^2$ exceeds $\eps>0$ at more
than ${c\over\eps\delta}$ times $t$ is bounded by $\delta$.
I sometimes loosely call this the number of errors.

\paradot{Sequence prediction}
Given a sequence $x_1x_2...x_{t-1}$,
we want to predict its likely continuation $x_t$. I assume that
the strings which have to be continued are drawn from a ``true''
probability distribution $\mu$.
The maximal prior information a prediction algorithm can possess
is the exact knowledge of $\mu$, but often the true distribution
is unknown. Instead, prediction is based on a guess $\rho$ of
$\mu$. While I require $\mu$ to be a measure, I allow $\rho$ to
be a semimeasure \cite{Li:97,Hutter:04uaibook}:\footnote{Readers
unfamiliar or uneasy with {\em semi}measures can without loss
ignore this technicality.}
Formally, $\rho:\X^*\to[0,1]$ is a semimeasure if
$\rho(x)\geq\sum_{a\in\X}\rho(xa)\,\forall x\in\X^*$, and a
(probability) measure if equality holds and $\rho(\epstr)=1$,
where $\epstr$ is the empty string. $\rho(x)$ denotes the
$\rho$-probability that a sequence starts with string $x$.
Further, $\rho(a|x):=\rho(xa)/\rho(x)$ is the ``posterior'' or ``predictive''
$\rho$-probability that the next symbol is $a\in\X$, given
sequence $x\in\X^*$.

\paradot{Bayes mixture}
We may know or assume that $\mu$ belongs to some countable class
${\cal M}:=\{\nu_1,\nu_2,...\}\ni\mu$ of semimeasures.
Then we can use the
weighted average on $\cal M$ (Bayes-mixture, data evidence, marginal)
\beq\label{xidef}
  \xi(x) :=
  \sum_{\nu\in\cal M}w_\nu\!\cdot\!\nu(x),\quad
  \sum_{\nu\in\cal M}w_\nu \leq 1,\quad w_\nu>0
\eeq
for prediction. One may interpret $w_\nu=\P[H_\nu]$ as prior belief
in $\nu$ and $\xi(x)=\P[x]$ as the subjective probability of $x$,
and $\mu(x)=\P[x|\mu]$ is the sampling distribution or likelihood.
The most important property of semimeasure $\xi$ is its dominance
\beq\label{xidom}
  \xi(x) \;\geq\; w_\nu\nu(x) \quad \forall x \;\mbox{and}\; \forall\nu\!\in\!\M,
  \qmbox{in particular} \xi(x) \;\geq\; w_\mu\mu(x)
\eeq
which is a strong form of absolute continuity.

\paradot{Convergence for deterministic environments}
In the predictive setting we are not interested in identifying the
true environment, but to predict the next symbol well. Let us
consider deterministic $\mu$ first. An environment is called
deterministic if $\mu(\a_{1:n})=1 \forall n$ for some sequence
$\a$, and $\mu=0$ elsewhere (off-sequence). In this case we
identify $\mu$ with $\a$ and the following holds:
\beq\label{detxibnd}
  \sum_{t=1}^\infty|1\!-\!\xi(\a_t|\a_{<t})| \;\leq\; \ln w_\a^{-1} \qqmbox{and}
  \xi(\a_{t:n}|\a_t)\to 1 \qmbox{for} n\geq t\to \infty
\eeq
where $w_\a>0$ is the weight of $\a\widehat=\mu\in\M$. This shows
that $\xi(\a_t|\a_{<t})$ rapidly converges to 1 and hence also
$\xi(\bar \a_t|\a_{<t})\to 0$ for $\bar\a_t\neq\a_t$, and that
$\xi$ is also a good multi-step lookahead predictor.
%
Proof: $\xi(\a_{1:n})\to c>0$, since $\xi(\a_{1:n})$ is monotone
decreasing in $n$ and $\xi(\a_{1:n})\geq
w_\mu\mu(\a_{1:n})=w_\mu>0$. Hence $\xi(\a_{1:n})/\xi(\a_{1:t})\to
c/c=1$ for any limit sequence $t,n\to\infty$. The bound follows from
$\sum_{t=1}^n 1-\xi(x_t|x_{<t}) \leq - \sum_{t=1}^n\ln
\xi(x_t|x_{<t}) = -\ln \xi(x_{1:n})$ and $\xi(\a_{1:n})\geq w_\a$.

\paradot{Convergence in probabilistic environments}
In the general probabilistic case we want to know how close
$\xi(x_t|x_{<t})$ is to the true probability $\mu(x_t|x_{<t})$.
One convenient distance measure is the (squared) Hellinger
distance
\beq\label{hdef}
  h_t(\o_{<t}) \;:=\; \sum_{a\in\X}(\sqrt{\xi(a|\o_{<t})}-\sqrt{\mu(a|\o_{<t})})^2
\eeq
One can show \cite{Hutter:03spupper,Hutter:04uaibook} that
\beq\label{hbnd}
  \sum_{t=1}^n\E{\textstyle\left[\!\left(\sqrt{{\xi(\o_t|\o_{<t})
       \over\mu(\o_t|\o_{<t})}}\!-\!1\right)^2\right]}
  \;\leq\; \sum_{t=1}^n\E[h_t]
  \;\leq\; D_n(\mu||\xi) := \textstyle \E[\ln{\mu(\o_{1:n})\over\xi(\o_{1:n})}]
  \;\leq\; \ln w_\mu^{-1}
\eeq
The first two inequalities actually hold for any two (semi)measures, and
the last inequality follows from \req{xidom}.
These bounds (with $n=\infty$) imply $h_t\to 0$ and hence
\ifjournal
\beqn
 \xi(x_t|\o_{<t})-\mu(x_t|\o_{<t}) \;\longrightarrow\; 0
 \qmbox{for any} x_t \quad\qmbox{and}\quad
 {\xi(\o_t|\o_{<t})\over\mu(\o_t|\o_{<t})} \;\longrightarrow\; 1
\eeqn
both rapid w.p.1 for $t\to\infty$.
\else
\beqn
 \mbox{$\xi(x_t|\o_{<t})-\mu(x_t|\o_{<t})\to 0$
 for any $x_t$ and
 ${\xi(\o_t|\o_{<t})\over\mu(\o_t|\o_{<t})} \to 1$, both
 rapid w.p.1 for $t\to\infty$}
\eeqn
\fi
An improved bound $\E[\exp(\odt\sum_t h_t)] \leq w_\mu^{-1/2}$
\cite{Hutter:04mlconvx} even shows that the probability that $\sum_t
h_t$ additively exceeds $\ln w_\mu^{-1}$ by $c$ (e.g.\ $c>10$) is
tiny $\e^{-c/2}$. One can also show multi-step lookahead convergence
$\xi(x_{t:n_t}|\o_{<t})-\mu(x_{t:n_t}|\o_{<t})\to 0$ (even for
unbounded horizon $1\leq n_t-t+1\to\infty$), which is interesting
for delayed sequence prediction and in reactive environments
\cite{Hutter:04uaibook}.
Since $\xi$ rapidly converges to $\mu$, one can anticipate that
also decisions based on $\xi$ are good.

\paradot{Bayesian decisions}
Let $\ell_{x_t y_t}\in[0,1]$ be the received loss when predicting
$y_t\in\cal Y$, but $x_t\in\cal X$ turns out to be the true
$t^{th}$ symbol of the sequence.
The $\rho$-optimal predictor
\beq\label{xlrdef}
  y_t^{\smash{\Lambda_\rho}}(\o_{<t}) \;:=\; \arg\min_{y_t}\sum_{x_t}\rho(x_t|\o_{<t})\ell_{x_t y_t}
\eeq
minimizes the $\rho$-expected loss.
For instance for ${\cal X}={\cal Y}=\{0,1\}$, $\Lambda_\rho$ is a
threshold strategy with $y_t^{\smash{\Lambda_\rho}}=0/1$
for $\rho(1|\o_{<t})\,_<^>\,\gamma$, where
$\gamma:={\ell_{01}-\ell_{00} \over
\ell_{01}-\ell_{00}+\ell_{10}-\ell_{11}}$.
The instantaneous loss at time
$t$ and the total $\mu$(=true)-expected loss for the first $n$
symbols are
\beq\label{ldef}
  l_t^{\smash{\Lambda_\rho}}(\o_{<t}) \;:=\; \E[\ell_{\o_t
  y_t^{\smash{\Lambda_\rho}}}|\o_{<t}] \qqmbox{and}
  L_n^{\smash{\Lambda_\rho}} \;:=\; \sum_{t=1}^n\E[\ell_{\o_t
  y_t^{\smash{\Lambda_\rho}}}]
\eeq
Let $\Lambda$ be {\em any} prediction scheme (deterministic or
probabilistic) with no constraint at all, taking {\em any} action
$y_t^\Lambda\in\cal Y$ with total expected loss $L_n^\Lambda$. If
$\mu$ is known, $\Lambda_\mu$ is obviously the best prediction
scheme in the sense of achieving minimal expected loss
$L_n^{\smash{\Lambda_\mu}}\leq L_n^\Lambda$ for any $\Lambda$. For
the predictor $\Lambda_\xi$ based on the Bayes mixture $\xi$, one
can show (proof in Appendix \ref{secApp}; see also
\cite{Merhav:98,Hutter:03spupper} for related bounds)
\beq\label{lbnd}
  (\sqrt{L_n^{\smash{\Lambda_\xi}}} - \sqrt{L_n^{\smash{\Lambda_\mu}}}\,)^2
  \;\leq\; \sum_{t=1}^n\E[(\sqrt{l_t^{\smash{\Lambda_\xi}}} - \sqrt{l_t^{\smash{\Lambda_\mu}}}\,)^2]
  \;\leq\; \sum_{t=1}^n 2\E[h_t]
\eeq
which actually holds for any two (semi)measures. Chaining with
\req{hbnd} implies, for instance, $l_t^{\smash{\Lambda_\xi}}\to
l_t^{\smash{\Lambda_\mu}}$ rapid w.p.1,
$\sqrt{L_n^{\smash{\Lambda_\xi}}}$ exceeds
$\sqrt{L_n^{\smash{\Lambda_\mu}}}$ by at most $\sqrt{2\ln
w_\mu^{-1}}$,
$L_n^{\smash{\Lambda_\xi}}/L_n^{\smash{\Lambda_\mu}}\to 1$ for
$L_n^{\smash{\Lambda_\mu}}\to\infty$, or if
$L_\infty^{\smash{\Lambda_\mu}}$ is finite, then also
$L_\infty^{\smash{\Lambda_\xi}}$. This shows that $\xi$ (via
$\Lambda_\xi$) performs also excellent from a decision-theoretic
perspective, i.e.\ suffers loss only slightly larger than the
optimal $\Lambda_\mu$ predictor.

One can also show that {\ir $\Lambda_\xi$ is Pareto-optimal}
(admissible) in the sense that every other predictor with smaller
loss than $\Lambda_\xi$ in some environment $\nu\in\M$ must
be worse in another environment \cite{Hutter:03optisp}.

\paradot{Continuous environmental classes}
I will argue later that countable $\cal M$ are sufficiently large
from a philosophical and computational perspective. On the other
hand, countable $\M$ exclude all continuously parameterized families
(like the class of all i.i.d. or Markov processes), common in
statistical practice. I show that the bounds above remain
approximately valid for most parametric model classes. Let
\beqn
  \M \;:=\; \{\nu_\th:\th\in\Theta\subseteq \SetR^d\}
\eeqn
be a family of probability distributions parameterized by a
$d$-dimensional continuous parameter $\th$, and $\mu \equiv
\nu_{\th_0} \in \M$ the true generating distribution. For a
continuous weight density\footnote{$w()$ will always denote
densities, and $w_{()}$ probabilities.} $w(\th)>0$
the sums \req{xidef} are
naturally replaced by integrals:
\beq\label{xidefc}
  \xi(x_{1:n}) \;:=\; \int_\Theta \!
  w(\th)\!\cdot\!\nu_\th(x_{1:n})\,d\th, \qquad
  \int_\Theta \! w(\th) \,d\th \;=\; 1 \qquad
\eeq
The most important property of $\xi$ was the dominance \req{xidom}
achieved by dropping the sum over $\nu$. The analogous
construction here is to restrict the integral over $\th$ to a
small vicinity of $\th_0$. Since a continuous parameter can
typically be estimated to accuracy $\propto n^{-1/2}$ after $n$
observations, the largest volume in which $\nu_\th$ as a function of $\th$ is
approximately flat is $\propto (n^{-1/2})^d$, hence
$\xi(x_{1:n})\gtrsim n^{-d/2}w(\th_0)\mu(x_{1:n})$. Under some
weak regularity conditions one can prove
\cite{Clarke:90,Hutter:03optisp}
\beq\label{cbnd}
  D_n(\mu||\xi) \;:=\;
  \E\textstyle
  \ln{\mu(\o_{1:n}) \over \xi(\o_{1:n})} \;\;\leq\;\;
  \ln w(\th_0)^{-1} + {d\over 2}\ln{n\over 2\pi} +
  {1\over 2}\ln\det\bar\jmath_n(\th_0) + o(1)
\eeq
where $w(\th_0)$ is the weight density (\ref{xidefc}) of $\mu$ in
$\xi$, and $o(1)$ tends to zero for $n \to\infty$, and the average
Fisher information matrix $\bar\jmath_n(\th) = -{1\over
n}\E[\nabla_\th\nabla^T_\th\ln\nu_\th(\o_{1:n})]$ measures the
local smoothness of $\nu_\th$ and is bounded for many reasonable
classes, including all stationary ($k^{th}$-order) finite-state
Markov processes. See Section \ref{secIID} for an application to
the i.i.d.\ ($k=0$) case.
We see that in the continuous case, $D_n$ is no longer bounded by
a constant, but grows very slowly (logarithmically) with $n$,
which still implies that $\eps$-deviations are exponentially
seldom.
Hence, \req{cbnd} allows to bound \req{hbnd} and \req{lbnd} even
in case of continuous $\M$.

\section{How to Choose the Prior}\label{secPrior}

I showed in the last section how to predict if the true
environment $\mu$ is unknown, but known to belong some class $\M$
of environments. In this section, I assume $\M$ to be given, and
discuss how to (universally) choose the prior $w_\nu$. After
reviewing various classical principles (indifference, symmetry,
minimax) for obtaining objective priors for ``small'' $\M$, I
discuss large $\M$. Occam's razor in conjunction with Epicurus'
principle of multiple explanations, quantified by Kolmogorov
complexity, leads us to a universal prior, which results in a
better predictor than any other prior over countable $\M$.

\paradot{Classical principles}
The probability axioms (implying Bayes' rule) allow to compute
posteriors and predictive distributions from prior ones, but are mute
about how to choose the prior. Much has been written on the choice of priors
(see \cite{Kass:96} for a survey and references).
A main classification is between objective and subjective priors.
An {\em objective prior} $w_\nu$ is a prior constructed based on some
rational principles, which ideally everyone without (relevant)
extra prior knowledge should adopt. In contrast, a {\em subjective
prior} aims at modelling the agents personal (subjective)
belief in environment $\nu$ prior to observation of $x$, but based
on his past personal experience or knowledge (e.g.\ of related
phenomena). In Section \ref{secDisc},
I show that one way to arrive at a subjective prior is to start
with an objective prior, make all past personal experience
explicit, determine a ``posterior'' and use it as subjective
prior. So I concentrate in the following on the more important
objective priors.

Consider a very simple case of two environments, e.g.\ a biased
coin with head {\em or} tail probability $1/3$. In absence of any
extra knowledge (which I henceforth assume) there is no reason to
prefer head probability $\th=1/3$ over $\th=2/3$ and vice versa, leaving
$w_{1/3}=w_{2/3}=\odt$ as the only rational choice.
More generally, for finite $\M$, the {\em symmetry or indifference
argument} \cite{Laplace:1812} suggests to set $w_\nu={1\over|\M|}$
$\forall\nu\in\M$. Similarly for a compact measurable parameter
space $\Theta$ we may choose a uniform density
$w(\th)=[\mbox{Vol}(\Theta)]^{-1}$. But there is a problem: If we
go to a different parametrization (e.g.\ $\th\leadsto\thp:=\sqrt{\th}$
in the Bernoulli case), the prior $w(\th)\leadsto w'(\thp)$
becomes non-uniform. Jeffreys' \cite{Jeffreys:46} solution is to
find a symmetry group of the problem (like permutations for finite
$\M$) and require the prior to be {\em invariant under group
transformations}. For instance, if $\th\in\SetR$ is a location
parameter (e.g.\ the mean) it is natural to require a
translation-invariant prior. Problems are that there may be no
obvious symmetry, the resulting prior may be improper (like
for the translation group), and the result can depend on which
parameters are treated as nuisance parameters.

The {\em maximum entropy principle} extends the symmetry principle by
allowing certain types of constraints on the parameters.
{\em Conjugate priors} are classes of priors such that the posteriors are
themselves again in the class. While this can lead to interesting
classes, the principle itself is not selective, since e.g.\ the
class of all priors forms a conjugate class.

Another {\em minimax approach} by Bernardo
\cite{Bernardo:79,Clarke:90} is to consider bound \req{cbnd},
which can actually be improved within $o(1)$ to an equality. Since
we want $D_n$ to be small, we minimize the r.h.s.\ for the worst
$\mu\in\M$. Choice $w(\th)\propto\sqrt{\det\bar\jmath_n(\th)}$
equalizes and hence minimizes \req{cbnd}. The problems are the
same as for Jeffrey's prior (actually often both priors coincide),
and also the dependence on the model class and potentially on $n$.

The principles above, although not unproblematic, {\em can}
provide good objective priors in many cases of small
discrete or compact spaces, but we will meet some more problems
later. For ``large'' model classes I am interested in, i.e.\
countably infinite, non-compact, or non-parametric spaces, the
principles typically do not apply or break down.

\paradot{Occam's razor et al}
Machine learning, the computer science branch of statistics, often
deals with very large model classes. Naturally, machine learning
has (re)discovered and exploited quite different principles for
choosing priors, appropriate for this situation. The overarching
principles put
together by Solomonoff \cite{Solomonoff:64} are: %
Occam's razor (choose the simplest model consistent with the data), %
Epicurus' principle of multiple explanations (keep all explanations consistent with the data), %
(Universal) Turing machines (to compute, quantify and assign codes to all quantities of interest), and %
Kolmogorov complexity (to define what simplicity/complexity means).

I will first ``derive'' the so called universal prior, and
subsequently justify it by presenting various welcome theoretical
properties and by examples.
The idea is that a priori, i.e.\ before seeing the data, all
models are ``consistent,$\!$'' so a-priori Epicurus would regard all
models (in $\M$) possible, i.e.\ choose $w_\nu>0$
$\forall\nu\in\M$. In order to also do (some) justice to Occam's
razor we should {\em prefer} simple hypotheses, i.e.\ assign high
prior (low) prior $w_\nu$ to simple (complex) hypotheses $H_\nu$. Before I
can define this prior, I need to quantify the notion of
complexity.

\paradot{Notation}
A function $f:\S\to\SetR\cup\{\pm\infty\}$ is said to be
lower semi-computable (or enumerable) if the set
$\{(x,y)\,:\,y<f(x),\, x\in\S,\, y\in\SetQ\}$ is recursively
enumerable. $f$ is upper semi-computable (or co-enumerable) if
$-f$ is enumerable. $f$ is computable (or recursive)
if $f$ and $-f$ are enumerable. The set of (co)enumerable
functions is recursively enumerable.
%
I write $O(1)$ for a constant of reasonable size: For instance, a
sequence of length $100$ is reasonable, maybe even $2^{30}$, but
$2^{500}$ is not. I write $f(x)\leqa g(x)$ for $f(x)\leq g(x)+O(1)$
and $f(x)\leqm g(x)$ for $f(x)\leq 2^{O(1)}\cdot g(x)$.
Corresponding equalities hold if the inequalities hold in both
directions.\footnote{I will ignore these additive and multiplicative
fudges in the discussion till Section \ref{secDisc}.}
%
We say that a property $A(n)\in\{true,false\}$ holds for {\em
most} $n$, if $\#\{t\leq n:A(t)\}/n\toinfty{n} 1$.

\paradot{Kolmogorov complexity}
We can now quantify the complexity of a string.
Intuitively, a string is simple if it can be described in a few
words, like ``the string of one million ones'', and is complex if
there is no such short description, like for a random object whose
shortest description is specifying it bit by bit. We are
interested in effective descriptions, and hence restrict decoders
to be Turing machines (TMs).
Let us choose some universal (so-called prefix) {\em Turing
machine $U$} with binary input=program tape, $\X$ary output tape, and
bidirectional work tape. We can then define the
{\em prefix Kolmogorov complexity}
\cite{Chaitin:75,Gacs:74,Kolmogorov:65,Levin:74} of
string $x$ as the length $\l$ of the shortest binary program $p$ for
which $U$ outputs $x$:
\beqn
  K(x) \;:=\; \min_p\{\l(p): U(p)=x\}
\eeqn
Simple strings like 000...0 can be generated by short programs,
and, hence have low Kolmogorov complexity, but irregular (e.g.\
random) strings are their own shortest description, and hence have
high Kolmogorov complexity.
For non-string objects $o$ (like numbers and functions) we define
$K(o):=K(\langle o\rangle)$, where $\langle o\rangle\in\X^*$ is
some standard code for $o$. In particular, if $(f_i)_{i=1}^\infty$ is
an enumeration of all (co)enumerable functions, we define
$K(f_i)=K(i)$.

An important property of $K$ is that it is nearly independent of
the choice of $U$. More precisely, if we switch from one universal
TM to another, $K(x)$ changes at most by an additive constant
independent of $x$. For natural universal TMs, the compiler
constant is of reasonable size $O(1)$.
A defining property of $K:\X^*\to\SetN$ is that it additively
dominates all co-enumerable functions $f:\X^*\to\SetN$ that
satisfy Kraft's inequality $\sum_x 2^{-f(x)}\;\leq\;1$, i.e.\
$K(x)\leqa f(x)$ for $K(f)=O(1)$. The universal TM provides a
shorter prefix code than any other effective prefix code.
$K$ shares many properties with Shannon's entropy
(information measure) $S$, but $K$ is superior to $S$ in many
respects. To be brief, $K$ is an excellent universal complexity
measure, suitable for quantifying Occam's razor.
We need the following properties of $K$:
\begin{list}{$\bullet$}{\parskip=0ex\parsep=0ex\itemsep=0ex\topsep=1ex}
\item[$a)$] $K$ is not computable, but only upper semi-computable,
\item[$b)$] the upper bound $K(n)\leqa \lb n+2\lb\log n$, \hfill (\refstepcounter{equation}\theequation\label{Kprop}) %
\item[$c)$] Kraft's inequality $\sum_x 2^{-K(x)}\leq1$, %
      which implies $2^{-K(n)}\leq\odn$ for most $n$, %
\item[$d)$] \ifjournal inform.\ \else information \fi non-increase $K(f(x))\leqa K(x)+K(f)$ for recursive $f:\X^*\to\X^*$, %
\item[$e)$] \ifjournal the MDL bound \fi $K(x)\leqa -\lb P(x)+K(P)$ \ifjournal \\ \fi
            if $P:\X^*\to[0,1]$ is enumerable and $\sum_x P(x)\leq 1$, %
\item[$f)$] $\sum_{x:f(x)=y} 2^{-K(x)}\eqm 2^{-K(y)}$ if $f$ is recursive and $K(f)=O(1)$.
\end{list}
The proof of $(f)$ can be found in Appendix \ref{secApp} and the
proofs of $(a)-(e)$ in \cite{Li:97}.

\paradot{The universal prior}
We can now quantify a prior biased towards simple models. First,
we quantify the complexity of an environment $\nu$ or hypothesis
$H_{\nu}$ by its Kolmogorov complexity $K(\nu)$. The universal
prior should be a decreasing function in the model's
complexity, and of course sum to (less than) one. Since
$K$ satisfies Kraft's inequality (\ref{Kprop}$c$), this suggests the following
choice:
\beq\label{uprior}
  w_\nu \;=\; w^U_\nu \;:=\; 2^{-K(\nu)}
\eeq
For this choice, the bound \req{hbnd} on $D_n$ (which bounds
\req{hbnd} and \req{lbnd}) reads
\beq\label{DnKCbnd}\textstyle
  \sum_{t=1}^\infty \E[h_t] \;\leq\; D_\infty \;\leq\; K(\mu)\ln 2
\eeq
i.e.\ the number of times, $\xi$ deviates from $\mu$ or
$l^{\Lambda_\xi}$ deviates from $l^{\Lambda_\mu}$ by more than
$\eps>0$ is bounded by $O(K(\mu))$, i.e.\ is proportional to the
complexity of the environment. Could other choices for $w_\nu$
lead to better bounds? The answer is essentially no
\cite{Hutter:04uaibook}: Consider any other reasonable prior
$w'_\nu$, where reasonable means (lower semi)computable with a
program of size $O(1)$. Then, MDL bound (\ref{Kprop}$e$) with
$P()\leadsto w'_{()}$ and $x\leadsto\langle\mu\rangle$ shows
$K(\mu)\leqa-\lb w'_\mu+K(w'_{()})$, hence $\ln w_\mu'\!^{-1}\geqa
K(\mu)\ln 2$ leads (within an additive constant) to a weaker
bound.
A counting argument also shows that $O(K(\mu))$ errors for most
$\mu$ are unavoidable. So this choice of prior leads to
very good prediction.

Even for continuous classes $\M$, we can assign a (proper) universal prior
(not density) $w_\th^U=2^{-K(\th)}>0$ for computable $\th$, and 0
for uncomputable ones. This effectively reduces $\M$ to a discrete
class $\{\nu_\th\in\M:w_\th^U>0\}$ which is typically dense in $\M$.
We will see that this prior has many advantages over the classical
prior densities.

\section{Independent Identically Distributed Data}\label{secIID}

I now compare the classical continuous prior densities to the
universal prior on classes of i.i.d.\ environments. I present some
standard critiques to the former, illustrated on Bayes-Laplace's
classical Bernoulli class with uniform prior: the problem of zero
p(oste)rior, non-confirmation of universal hypotheses, and
reparametrization and regrouping non-invariance. I show that the
universal prior does not suffer from these problems. Finally I
complement the general total bounds of Section \ref{secBSP} with
some i.i.d.-specific instantaneous bounds.

\paradot{Laplace's rule for Bernoulli sequences}
Let $x=x_1x_2...x_n\in\X^n=\B^n$ be generated by a biased coin
with head=1 probability $\th\in[0,1]$, i.e.\ the likelihood of $x$
under hypothesis $H_\th$ is
$\nu_\th(x)=\P[x|H_\th]=\th^{n_1}(1-\th)^{n_0}$, where
$n_1=x_1+...+x_n=n-n_0$. Bayes \cite{Bayes:1763} assumed a uniform
prior density $w(\th)=1$. The evidence is $\xi(x)=\int_0^1
\nu_\th(x)w(\th)\,d\th={n_1!n_0!\over(n+1)!}$ and the posterior
probability weight density $w(\th|x)=\nu_\th(x)w(\th)/\xi(x)
={(n+1)!\over n_1!n_0!}\th^{n_1}(1-\th)^{n_0}$ of $\th$ after
seeing $x$ is strongly peaked around the frequency estimate
$\hat\th={n_1\over n}$ for large $n$. Laplace \cite{Laplace:1812}
asked for the predictive probability $\xi(1|x)$ of observing
$x_{n+1}=1$ after having seen $x=x_1...x_n$, which is
$\xi(1|x)={\xi(x1)\over \xi(x)}={n_1+1\over n+2}$. (Laplace
believed that the sun had risen for $5\,000$ years = $1\,826\,213$ days
since creation, so he concluded that the probability of doom,
i.e.\ that the sun won't rise tomorrow is ${1\over 1826215}$.)
This looks like a reasonable estimate, since it is close to the
relative frequency, asymptotically consistent, symmetric, even
defined for $n=0$, and not overconfident (never assigns
probability 1).

\paradot{The problem of zero prior}
But also Laplace's rule is not without problems. The
appropriateness of the uniform prior has been questioned in
Section \ref{secPrior} and will be detailed below. Here I discuss
a version of the zero prior problem. If the prior is zero, then
the posterior is necessarily also zero. The above example seems
unproblematic, since the prior and posterior {\em densities}
$w(\th)$ and $w(\th|x)$ are non-zero. Nevertheless it is
problematic e.g.\ in the context of scientific confirmation theory
\cite{Earman:93}.

Consider the hypothesis $H$ that all balls in some urn, or all
ravens, are black (=1). A natural model is to assume that
balls (or ravens) are drawn randomly from an infinite population with
fraction $\th$ of black balls (or ravens) and to assume a uniform prior
over $\th$, i.e.\ just the Bayes-Laplace model. Now we draw $n$
objects and observe that they are all black.

We may formalize $H$ as the hypothesis $H':=\{\th=1\}$. Although
the posterior probability of the relaxed hypothesis
$H_\eps:=\{\th\geq 1-\eps\}$, $\P[H_\eps|1^n]=\int_{1-\eps}^1
w(\th|1^n)\,d\th=\int_{1-\eps}^1(n+1)\th^n d\th=1-(1-\eps)^{n+1}$
tends to 1 for $n\to\infty$ for every fixed $\eps>0$,
$\P[H'|1^n]=\P[H_0|1^n]$ remains identically zero, i.e.\ no amount
of evidence can confirm $H'$. The reason is simply that zero prior
$\P[H']=0$ implies zero posterior.

Note that $H'$ refers to the unobservable quantity $\th$ and only
demands blackness with probability 1. So maybe a better
formalization of $H$ is purely in terms of observational
quantities: $H'':=\{\o_{1:\infty}=1^\infty\}$. Since
$\xi(1^n)={1\over n+1}$, the predictive probability of observing
$k$ further black objects is
$\xi(1^k|1^n)={\smash{\xi(1^{n+k})\over\xi(1^n)}}={n+1\over n+k+1}$. While
for fixed $k$ this tends to 1,
$\P[H''|1^n]=\lim_{k\to\infty}\xi(1^k|1^n)\equiv 0$ $\forall n$,
as for $H'$.

One may speculate that the crux is the infinite population. But
for a finite population of size $N$ and sampling with (similarly
without) repetition, $\P[H''|1^n]=\xi(1^{N-n}|1^n)={n+1\over N+1}$
is close to one only if a large fraction of objects has been
observed. This contradicts scientific practice: Although only a
tiny fraction of all existing ravens have been observed, we regard
this as sufficient evidence for believing strongly in $H$.
This quantifies \cite[Thm.11]{Maher:04} and shows that Maher does
{\em not} solve the problem of confirmation of universal hypotheses.

There are two solutions of this problem: We may abandon
strict\eqbr/logical\eqbr/all-quantified\eqbr/universal hypotheses
altogether in favor of soft hypotheses like $H_\eps$. Although not
unreasonable, this approach is unattractive for several reasons.
The other solution is to assign a non-zero prior to $\th=1$.
Consider, for instance, the improper density
$w(\th)=\odt[1+\delta(1-\th)]$, where $\delta$ is the Dirac-delta
($\int f(\th)\delta(\th-a)\,d\th=f(a)$), or equivalently
$\P[\th\geq a]=1-\odt a$. We get
$\xi(x_{1:n})=\odt[{n_1!n_0!\over(n+1)!}+\delta_{0n_0}]$, where
$\delta_{ij}=\{{1 \text{ if } i=j\atop 0 \text{ else}\;\;\; }\}$
is Kronecker's $\delta$. In particular $\xi(1^n)=\odt{n+2\over
n+1}$ is much larger than for uniform prior. Since
$\xi(1^k|1^n)={n+k+2\over n+k+1}\cdot{n+1\over n+2}$, we get
$\P[H''|1^n]=\lim_{k\to\infty}\xi(1^k|1^n)={n+1\over n+2}\to 1$,
i.e.\ $H''$ gets strongly confirmed by observing a reasonable
number of black objects. This correct asymptotics also follows
from the general result \req{detxibnd}. Confirmation of $H''$ is
also reflected in the fact that $\xi(0|1^n)={1\over(n+2)^2}$ tends
much faster to zero than for uniform prior, i.e.\ the confidence that
the next object is black is much higher. The power actually
depends on the shape of $w(\th)$ around $\th=1$. Similarly $H'$
gets confirmed: $\P[H'|1^n]=\mu_1(1^n)\P[\th=1]/\xi(1^n)={n+1\over
n+2}\to 1$.
On the other hand, if a single (or more) 0 are observed ($n_0>0$),
then the predictive distribution $\xi(\cdot|x)$ and posterior
$w(\th|x)$ are the same as for uniform prior.

The findings above remain qualitatively valid for i.i.d.\
processes over finite non-binary alphabet $|\X|>2$ and for
non-uniform prior.

Surely to get a generally working setup, we should also assign a
non-zero prior to $\th=0$ and to all other ``special'' $\th$, like
$\odt$ and ${1\over 6}$, which may naturally appear in a hypothesis,
like ``is the coin or die fair''. The natural continuation of this
thought is to assign non-zero prior to all computable $\th$. This is
another motivation for the universal prior $w_\th^U=2^{-K(\th)}$
\req{uprior} constructed in Section \ref{secPrior}. It is difficult
but not impossible to operate with such a prior
\cite{Hutter:04mdlspeed,Hutter:06mdlspeedx}. One may want to mix the
discrete prior $w_\nu^U$ with a continuous (e.g.\ uniform) prior
density, so that the set of non-computable $\th$ keeps a non-zero
density. Although possible, we will see that this is actually not
necessary.

\paradot{Reparametrization invariance}
Naively, the uniform prior is justified by the indifference
principle, but as discussed in Section \ref{secPrior}, uniformity
is not reparametrization invariant. For instance if in our
Bernoulli example we introduce a new parametrization
$\thp=\sqrt\th$, then the $\thp$-density $w'(\thp)=2\sqrt{\th}w(\th)$
is no longer uniform if $w(\th)=1$ is uniform.

More generally, assume we have some principle which leads to some
prior $w(\th)$. Now we apply the principle to a different
parametrization $\thp\in\Thp$ and get prior $w'(\th')$. Assume
that $\th$ and $\thp$ are related via bijection
$\th=f(\thp)$. Another way to get a $\thp$-prior is to
transform the $\th$-prior $w(\th)\leadsto \tilde w(\thp)$.
The reparametrization invariance principle (RIP) states that
$w'$ should be equal to $\tilde w$.

For discrete $\Theta$, simply $\tilde w_\thb=w_{f(\thp)}$, and a
uniform prior remains uniform ($w'_\thb=\tilde
w_\thb=w_\th={1\over|\Theta|}$) in any parametrization, i.e.\ the
indifference principle satisfies RIP in finite model classes.

In case of densities, we have $\tilde
w(\thp)=w(f(\thp)){df(\thp)\over d\thp}$, and the indifference
principle violates RIP for non-linear transformations $f$. But
Jeffrey's and Bernardo's principle satisfy RIP. For instance, in
the Bernoulli case we have $\bar\jmath_n(\th)={1\over\th}+{1\over
1-\th}$, hence $w(\th)={1\over\pi}[\th(1-\th)]^{-1/2}$ and
$w'(\thp)={1\over\pi}[f(\thp)(1-f(\thp))]^{-1/2}{df(\thp)\over
d\thp}=\tilde w(\thp)$.

Does the universal prior $w_\th^U=2^{-K(\th)}$ satisfy RIP? If we
apply the ``universality principle'' to a $\thp$-parametrization,
we get $w_\thb'\!\!^U=2^{-K(\thp)}$. On the other hand, $w_\th$
simply transforms to $\tilde w_\thb^U = w_{f(\thp)}^U =
2^{-K(f(\thp))}$ ($w_\th$ is a discrete (non-density) prior, which
is non-zero on a discrete subset of $\M$). For computable $f$ we
have $K(f(\thp))\leqa K(\thp)+K(f)$ by (\ref{Kprop}$d$), and
similarly $K(f^{-1}(\th))\leqa K(\th)+K(f)$ if $f$ is invertible.
Hence for simple bijections $f$ i.e.\ for $K(f)=O(1)$, we have
$K(f(\thp))\equa K(\thp)$, which implies $w_\thb'\!\!^U \eqm
\tilde w_\thb^U$, i.e.\ {\ir the universal prior satisfies RIP}
w.r.t.\ simple transformations $f$ (within a multiplicative
constant).

\paradot{Regrouping invariance}
There are important transformations $f$ which are {\em not}
bijections, which we consider in the following. A simple
non-bijection is $\th=f(\thp)=\thp^2$ if we consider
$\thp\in[-1,1]$. More interesting is the following example: Assume
we had decided not to record blackness versus non-blackness of
objects, but their ``color''. For simplicity of exposition assume
we record only whether an object is black or white or colored, i.e.\
$\X'=\{B,W,C\}$. In analogy to the binary case we use the
indifference principle to assign a uniform prior on
$\v\thp\in\Thp:=\Delta_3$, where
$\Delta_d:=\{\v\thp\in[0,1]^d:\sum_{i=1}^d\thp_i=1\}$, and
$\nu_\thb(x'_{1:n})=\prod_i{\thp_i}^{n_i}$. All inferences regarding
blackness (predictive and posterior) are identical to the binomial
model $\nu_\th(x_{1:n})=\th^{n_1}(1-\th)^{n_0}$ with $x'_t=B$
$\leadsto$ $x_t=1$ and $x'_t=W\,$or$\,C$ $\leadsto$ $x_t=0$ and
$\th=f(\v\thp)=\thp_B$ and
$w(\th)=\int_{\Delta_3}w'(\v\thp)\delta(\thp_B-\th)d\v\thp$.
Unfortunately, for uniform prior $w'(\v\thp)\propto 1$, $w(\th)\propto
1-\th$ is {\em not} uniform, i.e.\ the indifference principle is
{\em not} invariant under splitting/grouping, or general
regrouping. Regrouping invariance is regarded as a very important
and desirable property \cite{Walley:96}.

I now consider general i.i.d.\ processes
$\nu_\th(x)=\prod_{i=1}^d\th_i^{n_i}$. Dirichlet priors
$w(\th)\propto\prod_{i=1}^d\th_i^{\a_i-1}$ form a natural
conjugate class ($w(\th|x)\propto\prod_{i=1}^d\th_i^{n_i+\a_i-1}$)
and are the default priors for multinomial (i.i.d.)\ processes over
finite alphabet $\X$ of size $d$. Note that $\xi(a|x)={n_a+\a_a\over
n+\a_1+...+\a_d}$ generalizes Laplace's rule and coincides with
Carnap's \cite{Carnap:52} confirmation function.
Symmetry demands $\a_1=...=\a_d$; for instance $\a\equiv 1$ for
uniform and $\a\equiv\odt$ for Bernard-Jeffrey's prior. Grouping
two ``colors'' $i$ and $j$ results in a Dirichlet prior with
$\a_{i\&j}=\a_i+\a_j$ for the group. The only way to respect
symmetry under all possible groupings is to set $\a\equiv 0$. This
is Haldane's improper prior, which results in unacceptably
overconfident predictions $\xi(1|1^n)=1$. Walley \cite{Walley:96}
solves the problem that there is no single acceptable prior
density by considering sets of priors.

I now show that the universal prior $w_\th^U=2^{-K(\th)}$ is
invariant under regrouping, and more generally under all simple
(computable with complexity O(1)) even non-bijective
transformations. Consider prior $w'_\thb$. If $\th=f(\thp)$ then
$w'_\thb$ transforms to $\tilde
w_\th=\sum_{\thp:f(\thp)=\th}w'_\thb$ (note that for
non-bijections there is more than one $w'_\thb$ consistent with
$\tilde w_\th$). In $\thp$-parametrization, the universal prior
reads $w_\thb'\!\!^U=2^{-K(\thp)}$.
Using (\ref{Kprop}$f$) with $x=\langle\thp\rangle$ and
$y=\langle\th\rangle$ we get
\beqn
  \tilde w_\th^U \;=\; \sum_{\thp:f(\thp)=\th}2^{-K(\thp)}
  \;\eqm\; 2^{-K(\th)} \;=\; w_\th^U
\eeqn
i.e.\ {\ir the universal prior is general transformation and hence
regrouping invariant} (within a multiplicative constant) w.r.t.\
simple computable transformations $f$.

Note that reparametrization and regrouping invariance hold
for arbitrary classes $\M$ and are not limited to the i.i.d.\ case.

\paradot{Instantaneous bounds}
The cumulative bounds \req{hbnd} and \req{cbnd} stay valid for
i.i.d.\ processes, but instantaneous bounds are now also possible.
For i.i.d.\ $\M$ with continuous, discrete, and universal prior,
respectively, one can show (in preparation; see
\cite{Krichevskiy:98,Hutter:04mdlspeed,Hutter:06mdlspeedx} for
related bounds)
\beq\label{iIIDbnd}
  \E[h_n] \;\leqm\; \odn\ln w({\th_0})^{-1} \qmbox{and}
  \E[h_n] \;\leqm\; \odn\ln w_{\th_0}^{-1} \;=\; \odn K(\th_0)\ln 2
\eeq
Note that, if summed up over $n$, they lead to weaker cumulative bounds.

\section{Universal Sequence Prediction}\label{secUSP}

Section \ref{secPrior} derived the universal prior and Section
\ref{secIID} discussed i.i.d.\ classes. What remains and will be
done in this section is to find a universal class of environments,
namely Solomonoff-Levin's class of all (lower semi)computable
(semi)measures. The resulting universal mixture is equivalent to
the output distribution of a universal Turing machine with uniform
input distribution. The universal prior avoids the problem of old
evidence and the universal class avoids the necessity of updating
$\M$. I discuss the general total bounds of Section \ref{secBSP}
for the specific universal mixture, and supplement them with some
weak instantaneous bounds. Finally, I show that the universal
mixture performs better than classical continuous mixtures, even
in uncomputable environments.

\paradot{Universal choice of $\cal M$}
The bounds of Section \ref{secBSP} apply if $\M$ contains the true
environment $\mu$. The larger $\M$ the less restrictive is this
assumption. The class of all computable distributions, although
only countable, is pretty large from a practical point of view.
(Finding a non-computable physical system would overturn the Church-Turing thesis.)
It is the largest class, relevant from a computational point of
view. Solomonoff \cite[Eq.(13)]{Solomonoff:64} defined and studied
the mixture over this class.

One problem is that this class is not enumerable, since
the class of computable functions $f:\X^*\to\SetR$ is not enumerable
(halting problem), nor is it decidable whether a function is a
measure. Hence $\xi$ is completely incomputable. Levin
\cite{Zvonkin:70} had the idea to ``slightly'' extend the class
and include also lower semi-computable semimeasures. One can show
that this class $\M_U=\{\nu_1,\nu_2,...\}$ is enumerable, hence
\beq\label{xiUdef}
  \xi_U(x) \;=\; \sum_{\nu\in\M_U} w_\nu^U \nu(x)
\eeq
is itself lower semi-computable, i.e.\ $\xi_U\in\M_U$, which is a
convenient property in itself. Note that since ${1\over n\log^2
n}\leqm w_{\nu_n}^U\leq\odn$ for most $n$ by (\ref{Kprop}$b$) and
(\ref{Kprop}$c$), most $\nu_n$ have prior approximately reciprocal
to their index $n$, as advocated by Jeffreys
\cite[p238]{Jeffreys:61} and Rissanen \cite{Rissanen:83uniprior}.

In some sense $\M_U$ is the largest
class of environments for which $\xi$ is in some sense computable
\cite{Hutter:03unipriors,Hutter:06unipriorx}, but see
\cite{Schmidhuber:02gtm} for even larger classes. Note that
including non-semi-computable $\nu$ would not affect $\xi_U$,
since $w_\nu^U=0$ on such environments.

\paradot{The problem of old evidence}
An important problem in Bayesian inference in general and
(Bayesian) confirmation theory \cite{Earman:93} in particular is
how to deal with `old evidence' or equivalently with `new
theories'. How shall a Bayesian treat the case when some evidence
$E\widehat=x$ (e.g.\ Mercury's perihelion advance) is known
well-before the correct hypothesis/theory/model $H\widehat=\mu$ (Einstein's
general relativity theory) is found? How shall $H$ be added to the
Bayesian machinery a posteriori? What is the prior of $H$? Should
it be the belief in $H$ in a hypothetical counterfactual world in
which $E$ is not known? Can old evidence $E$ confirm $H$? After
all, $H$ could simply be constructed/biased/fitted towards
``explaining'' $E$.

The universal class $\M_U$ and universal prior $w_\nu^U$ formally
solve this problem: The universal prior of $H$ is $2^{-K(H)}$.
This is independent of $\M$ and of whether $E$ is known or not. If
we use $E$ to construct $H$ or fit $H$ to explain $E$, this will
lead to a theory which is more complex ($K(H)\geqa K(E)$) than a
theory from scratch ($K(H)=O(1)$), so cheats are automatically
penalized. There is no problem of adding hypotheses to $\M$ a
posteriori. Priors of old hypotheses are not affected. Finally,
$\M_U$ includes {\em all} hypotheses (including yet unknown or
unnamed ones) a priori. So at least theoretically, updating $\M$
is unnecessary.

\paradot{Other representations of $\xi_U$}
Definition \req{xiUdef} is somewhat complex, relying on
enumeration of semimeasures and Kolmogorov complexity. I now
approach $\xi_U$ from a different perspective.
Assume that our world is governed by a computable {\em
deterministic} process describable in $\leq l$ bits. Consider a
standard (not prefix) Turing machine $U'$ and programs $p$
generating environments starting with $x$. Let us pad all programs
so that they have length exactly $l$. Among the $2^l$ programs of
length $l$ there are $N_l(x):=\#\{p\in\{0,1\}^l:U'(p)=x*\}$ programs
consistent with observation $x$. If we regard all environmental
descriptions $p\in\{0,1\}^l$ a priori as equally  likely (Epicurus)
we should adopt the relative frequency $N_l(x)/2^l$ as our prior
belief in $x$. Since we do not know $l$ and we can pad every $p$
arbitrarily, we could take the limit
$M(x):=\lim_{l\to\infty}N_l(x)/2^l$ (which exists, since
$N_l(x)/2^l$ increases).
Or equivalently: $M(x)$ is the probability that $U'$ outputs a
string starting with $x$ when provided with uniform random noise on
the program tape.
Note that a uniform distribution is also used in the No Free Lunch
theorems \cite{Wolpert:97} to prove the impossibility of universal
learners, but in our case the uniform distribution is piped
through a universal Turing machine which defeats these negative
implications.
Yet another representation of $M$ is as follows: For every $q$ printing
$x*$ there exists a shortest prefix (called minimal) $p$ of $q$ printing
$x$. $p$ possesses $2^{l-\l(p)}$ prolongations to length $l$, all
printing $x*$. Hence all prolongations of $p$ together yield a
contribution $2^{l-\l(p)}/2^l=2^{-\l(p)}$ to $M(x)$. Let $U(p)=x*$
iff $p$ is a minimal program printing a string starting with $x$.
Then
\beq\label{Mdef}
  M(x) \;=\; \sum_{p:U(p)=x*} 2^{-\l(p)}
\eeq
which may be regarded as a $2^{-\l(p)}$-weighted mixture over all
computable deterministic environments $\nu_p$ ($\nu_p(x)=1$ if $U(p)=x*$ and
0 else).
Now, as a positive surprise, {\ir $M(x)$ coincides with
$\xi_U(x)$} within an irrelevant multiplicative constant. So it is
actually sufficient to consider the class of {\em deterministic}
semimeasures. The reason is that the probabilistic semimeasures
are in the convex hull of the deterministic ones, and so need not
be taken extra into account in the mixture.
One can also get an explicit enumeration of all lower
semi-computable semimeasures $\M_U=\{\nu_1,\nu_2,...\}$ by means of
$\nu_i(x):=\sum_{p:T_i(p)=x*}2^{-\l(p)}$, where
$T_i(p) \equiv U(\langle i\rangle p)$, $i=1,2,...$ is an enumeration of all
monotone Turing machines.

\paradot{Bounds for computable environments}
The bound \req{DnKCbnd} surely is applicable for $\xi=\xi_U$ and
now holds for {\em any} computable measure $\mu$. Within an
additive constant the bound is also valid for $M\eqm\xi$. That is,
{\ir $\xi_U$ and $M$ are excellent predictors with the only condition
that the sequence is drawn from any computable probability
distribution}. Bound \req{DnKCbnd} shows that the total number of
prediction errors is small.
Similarly to \req{detxibnd} one can show that
$\sum_{t=1}^n|1-M(x_t|x_{<t})| \leq \Km(x_{1:n})\ln 2$, where the
monotone complexity $\Km(x):=\min\{\l(p):U(p)=x*\}$ is defined as
the length of the shortest (nonhalting) program computing a string
starting with $x$ \cite{Zvonkin:70,Li:97,Hutter:04uaibook}.

If $x_{1:\infty}$ is a computable sequence, then
$\Km(x_{1:\infty})$ is finite, which implies {\ir $M(x_t|x_{<t})\to 1$
on every computable sequence}. This means that if the environment is
a computable sequence (whichsoever, e.g.\ $1^\infty$ or the digits
of $\pi$ or $\e$), after having seen the first few digits, $M$
correctly predicts the next digit with high probability, i.e.\ it
recognizes the structure of the sequence. In particular, observing
an increasing number of black balls or black ravens or sunrises,
$M(1|1^n)\to 1$ ($\Km(1^\infty)=O(1)$) becomes rapidly confident
that future balls and ravens are black and that the sun will rise
tomorrow.

Total bounds \req{detxibnd} and \req{DnKCbnd} are suitable in an
online setting, but {\em given} a fixed number of $n$ observations,
they give no guarantee on the next instance.

\paradot{Weak instantaneous bounds}
In Section \ref{secIID}, I derived good instantaneous bounds for
i.i.d.\ classes. For coin or die flips or balls drawn from an urn
this model is appropriate. But ornithologists do not really sample
ravens independently at random. Although not strictly valid, the
i.i.d.\ model may in this case still serve as a useful proxy for the
true process. But to model the rise of the sun as an i.i.d.\ process
is more than questionable. On the other hand it is plausible that
these examples (and other processes like weather or stock-market)
are governed by {\em some} (probabilistic) computable process. So
model class $\M_U$ and predictor $M$ seem appropriate. While
excellent total bounds \req{detxibnd} and \req{DnKCbnd} exist, the
essentially only {\ir instantaneous bound} I was able to derive
(proof in Appendix \ref{secApp}) is
\beq\label{iMbnd}
  2^{-K(n)} \;\leqm\; M(\bar x_n|x_{<n})
  \;\leqm\; 2^{2\Km(x_{1:n})-K(n)}
\eeq
valid for all $n$ and $x_{1:n}$ and $\bar x_n\neq x_n$. I discuss
the bound for the sequence $x_{1:\infty}=1^\infty$, but most of what
I say remains valid for any other computable sequence. Since
$\Km(1^n)=O(1)$, we get
\beqn
   M(0|1^n) \;\eqm\; 2^{-K(n)}
\eeqn
Since $2^{-K(n)}\leq\odn$ for most $n$, this shows that $M$ quickly
disbelieves in non-black objects and doomsday, similarly as in the
i.i.d.\ model, but now only for most $n$.

\paradot{Magic numbers}
This `most' qualification has interesting consequences:
$M(0|1^n)$ spikes up for simple $n$. So $M$ is cautious at magic
instance numbers, e.g.\ fears doom on day $2^{20}$ more than on a
comparable random day. While this looks odd and pours water on the
mills of prophets, it is not completely absurd. For instance,
major software problems have been anticipated for the magic date,
1st of January 2000. There are many other occasions, where
something happens at ``magic'' dates or instances; for instance
solar eclipses.

Also, certain processes in nature follow fast growing sequences like
those of the powers of two (e.g.\ the number of cells in an early
human embryo) or the Fibonacci numbers (e.g.\ the number of petals
or the arrangement of seeds in some flowers). Finally, that numbers
with low (Kolmogorov) complexity cause high probability in real data
bases can readily be verified by counting their frequency in the
world wide web with Google \cite{Cilibrasi:06}.

On the other hand, (returning to sequence prediction) on most simple
dates, nothing exceptional happens. Due to the total bound
$\sum_{n=0}^\infty M(0|1^n)\leq O(1)$, $M$ cannot spike up too much
too often. $M$ tells us to be more prepared but not to expect the
unexpected on those days. Another issue is that often we do not know
the exact start of the sequence. How many ravens exactly have
ornithologists observed, and how many days exactly did the sun rise
so far? In absence of this knowledge we need to Bayes-average over
the sequence length which will wash out the spikes.

\paradot{Universal is better than continuous $\M$}
Although I argued that incomputable environments $\mu$ can safely be
ignored, one may be nevertheless uneasy using Solomonoff's
$M\eqm\xi_U$ (\ref{Mdef}) if outperformed by a continuous mixture
$\xi$ (\ref{xidefc}) on such $\mu\in\M\setminus\M_U$, for instance
if $M$ would fail to predict a Bernoulli($\th$) sequence for
incomputable $\th$. Luckily this is not the case: Although
$\nu_\th()$ and $w_\th$ can be incomputable, the studied classes
$\M$ themselves, i.e.\ the two-argument function $\nu_{()}()$, and
the weight function $w_{()}$, and hence $\xi()$, are typically
computable (the integral can be approximated to arbitrary
precision). Hence $M(x)\eqm\xi_U(x)\geq 2^{-K(\xi)}\xi(x)$ by
\req{xiUdef} and $K(\xi)$ is often quite small. This implies for
{\em all} $\mu$
\beqn\textstyle
  D_n(\mu||M)
    \;\equiv\; \E[\ln\!{\mu(\o_{1:n})\over M(\o_{1:n})}]
    \;=\; \E[\ln\!{\mu(\o_{1:n})\over\xi(\o_{1:n})}]
    \!+\! \E[\ln\!{\xi(\o_{1:n})\over M(\o_{1:n})}]
    \;\leqa\; D_n(\mu||\xi) \!+\! K(\xi)\ln 2
\eeqn
So any bound \req{cbnd} for $D_n(\mu||\xi)$ is directly valid also
for $D_n(\mu||M)$, save an additive constant. That is, $M$ is
superior (or equal) to all computable mixture predictors $\xi$ based
on any (continuous or discrete) model class $\M$ and weight
$w(\th)$, even if environment $\mu$ is {\em not} computable.
Furthermore, while for essentially all parametric classes,
$D_n(\mu||\xi)\sim{d\over 2}\ln n$ grows logarithmically in $n$ for all
(incl.\ computable) $\mu\in\M$, $D_n(\mu||M)\leq K(\mu)\ln 2$ is
finite for computable $\mu$.
Bernardo's prior even implies a bound for $M$ that is uniform
(minimax) in $\th\in\Theta$. Many other priors based on reasonable
principles are argued for (see Section \ref{secPrior} and
\cite{Kass:96}). The above shows that $M$ is superior to all of
them.
Actually the existence of {\em any} computable probabilistic
predictor $\rho$ with $D_n(\mu||\rho)=o(n)$ is sufficient for
$M$ to predict $\mu$ equally well.

\paradot{Future bounds}
Another important question is how many errors are still to come
after some grace or learning period. Formally, given $x_{1:n}$,
how large is the future expected error
$r_n:=\sum_{t=n+1}^\infty\E[h_t|\o_{1:n}=x_{1:n}]$? The total
bound \req{hbnd}+\req{DnKCbnd} only implies that $r_n$
asymptotically tends to zero w.p.1, and the instantaneous
bounds \req{iIIDbnd} and \req{iMbnd} are weak and do not sum up finitely.
Since the complexity of $\mu$ bounds the total loss, a natural guess
is that something like the conditional complexity of $\mu$ given
$x$ (on an extra input tape) bounds the future loss.
Indeed one can show \cite{Hutter:04uaibook,Hutter:05postbnd}
\beq\label{pbnd}
  \sum_{t=n+1}^\infty \E[h_t|\o_{1:n}]
  \;\leqa\; [K(\mu|\o_{1:n})\!+\!K(n)]\ln 2
\eeq
i.e.\ if our past observations $\o_{1:n}$ contain a lot of
information about $\mu$, we make few errors in future.
For instance, consider the large space $\X$ of
pixel images, and all observations are identical
$\mu\,\widehat=\,\o=x_1x_1x_1...$, where $x_1$ is a ``typical'' image
of complexity, say, $K(x_1)\equa 10^6\equa\Km(\o)$. Obviously, after
seeing a couple of identical images we expect the next one to be
the same again. While total bound \req{DnKCbnd} quite uselessly tells
us that $M$ makes less than $10^6$ errors, future bound
\req{pbnd} with $n=1$ shows that $M$ makes only $K(\mu|x_1)=O(1)$
errors.
The $K(n)$ term can be improved to the complexity of the randomness
deficiency of $\o_{1:n}$ if a more suitable variant of algorithmic
complexity that is monotone in its condition is used
\cite{Hutter:05postbnd,Hutter:07postbndx}.
No future bounds analogous to \req{pbnd} for general prior or class
are known.

\section{Discussion}\label{secDisc}

\paradot{Critique and problems}
In practice we often have extra information about the problem at
hand, which could and should be used to guide the forecasting. One
way is to explicate all our prior knowledge $y$ and place it on an
extra input tape of our universal Turing machine $U$, which leads to
the conditional complexity $K(\cdot|y)$. We now assign
``subjective'' prior $w_{\nu|y}^U=2^{-K(\nu|y)}$ to environment
$\nu$, which is large for those $\nu$ that are simple (have short
description) relative to our background knowledge $y$. Since
$K(\mu|y)\leqa K(\mu)$, extra knowledge never misguides (see
\req{DnKCbnd}). Alternatively we could prefix our observation
sequence $x$ by $y$ and use $M(yx)$ for prediction
\cite{Hutter:04uaibook}.

Another critique concerns the dependence of $K$ and $M$ on $U$.
Predictions for short sequences $x$ (shorter than typical compiler
lengths) can be arbitrary. But taking into account our (whole)
scientific prior knowledge $y$, and predicting the now long string
$yx$ leads to good (less sensitive to ``reasonable'' $U$)
predictions \cite{Hutter:04uaibook}. For an interesting attempt to
make $M$ unique see \cite{Mueller:06}.

Finally, $K$ and $M$ can serve as ``gold standards'' which
practitioners should aim at, but since they are only
semi-computable, they have to be (crudely) approximated in
practice. Levin complexity \cite{Li:97}, the speed prior
\cite{Schmidhuber:02speed}, the minimal message and description
length principles \cite{Rissanen:89,Wallace:05}, and off-the-shelf
compressors like Lempel-Ziv \cite{Lempel:76} are such
approximations, which have been successfully applied to a plethora
of problems \cite{Cilibrasi:05,Schmidhuber:04oops}.

\paradot{Summary}
I compared traditional Bayesian sequence prediction based on
continuous classes and prior densities to Solomonoff's universal
predictor $M$, prior $w_\nu^U$, and class $\M_U$. I discussed the
following advantages (+) and problems ($-$) of Solomonoff's approach:
\vspace{-1ex}
\begin{itemize}\parskip=0ex\parsep=0ex\itemsep=0ex
\item[+] general total bounds for generic class, prior, and loss, %
\item[+] universal and i.i.d.-specific instantaneous and future bounds, %
\item[+] the $D_n$ bound for continuous classes, %
\item[+] indifference/symmetry principles, %
\item[+] the problem of zero p(oste)rior and confirmation of universal hypotheses, %
\item[+] reparametrization and regrouping invariance, %
\item[+] the problem of old evidence and updating, %
\item[+] that $M$ works even in non-computable environments, %
\item[+] how to incorporate prior knowledge,
\item[$-$] the prediction of short sequences, %
\item[$-$] the constant fudges in all results and the $U$-dependence, %
\item[$-$] $M$'s incomputability and crude practical approximations. %
\end{itemize}\vspace{-1ex}
In short, universal prediction solves or avoids or meliorates many
foundational and philosophical problems, but has to be compromised
in practice.

\paradot{Conclusion}
The goal of the paper was to establish a single, universal theory
for (sequence) prediction and (hypothesis) confirmation, applicable
to all inductive inference problems. I started by showing that
Bayesian prediction is consistent for any countable model class,
provided it contains the true distribution. The major (agonizing)
problem Bayesian statistics leaves open is how to choose the model
class and the prior. Solomonoff's theory fills this gap by choosing
the class of all computable (stochastic) models, and a universal
prior inspired by Ockham and Epicurus, and quantified by Kolmogorov
complexity. I discussed in breadth how and in which sense this
theory solves the inductive inference problem, by studying a
plethora of problems other approaches suffer from. In one line: All
you need for universal prediction is Ockham, Epicurus, Bayes,
Solomonoff, Kolmogorov, and Turing. By including Bellman, one can
extend this theory to universal decisions in reactive environments
\cite{Hutter:04uaibook}.

\paradot{Acknowledgements}
I would like to thank Frank Stephan for his detailed feedback
on earlier drafts.


\begin{small}

\end{small}

\pagebreak[3]\appendix
\section{Proofs of \req{lbnd}, (\ref{Kprop}$f$), and \req{iMbnd}}\label{secApp}

\paradot{Proof of loss bound \req{lbnd}}
Let $X$ and $Y$ be real-valued random variables.
Taking the square root of the well-known Schwarz inequality
$(\E[XY])^2\leq\E[X^2]\E[Y^2]$ we get
\beqn
  \E[(X\!-\!Y)^2] - (\sqrt{\E[X^2]}\!-\!\sqrt{\E[Y^2]}\,)^2
  \;\equiv\; 2\sqrt{\E[X^2]\E[Y^2]} - 2\E[XY] \;\geq\; 0.
\eeqn
Substituting $X\leadsto\sqrt{a_i}$, $Y\leadsto\sqrt{b_i}$,
$\E[...]\leadsto{1\over v_\Sigma}\sum_i v_i...$ with
$v_\Sigma:=\sum_i v_i$, we get, after multiplying with $v_\Sigma$, the
``Hellinger'' bound
\beq\label{eqHellBnd}\textstyle
  (\sqrt{\sum_i v_i a_i} - \sqrt{\sum_i v_i b_i\!}\,)^2
  \;\leq\;
  \sum_i v_i (\sqrt{a_i}-\sqrt{b_i\!}\,)^2
\eeq
for real $a_i,b_i,v_i\geq 0$ (also valid for $v_\Sigma=0$). I will
use \req{eqHellBnd} three times in proving \req{lbnd}.
With the abbreviations $m=y_t^{\smash{\Lambda_\mu}}$ and
$s=y_t^{\smash{\Lambda_\xi}}$ and
\beqn
  \X=\{1,...,N\},\quad
  N=|\X|, \quad
  i=x_t, \quad
  y_i=\mu(x_t|\o_{<t}), \quad
  z_i=\xi(x_t|\o_{<t})
\eeqn
the loss \req{ldef} and Hellinger distance \req{hdef} can then be
expressed by $l_t^{\smash{\Lambda_\xi}}=\sum_i y_i \ell_{is}$,
$l_t^{\smash{\Lambda_\mu}}=\sum_i y_i \ell_{im}$ and $h_t=\sum_i
(\sqrt{z_i}-\sqrt{y_i})^2$. By definition \req{xlrdef} of
$y_t^{\smash{\Lambda_\mu}}$ and $y_t^{\smash{\Lambda_\xi}}$ we have
\beq\label{lcnstr2}
 \sum_i y_i \ell_{im}\leq\sum_i y_i \ell_{ij} \qmbox{and}
 \sum_i z_i \ell_{is}\leq\sum_i z_i \ell_{ij} \qmbox{for all} j.
\eeq
Actually, I need the first constraint only for $j=s$ and the
second for $j=m$.
From (\ref{lcnstr2}) we get
\bqa\label{eqsqrtl}
  & & \nq\textstyle \sqrt{\sum_i y_i\l_{is}} - \sqrt{\sum_i y_i\l_{im}} \;\geq\; 0 \qqmbox{and}
\\ \nonumber
  & & \nq \textstyle[{\partial\over\partial\l_{is}}\!+\!{\partial\over\partial\l_{im}}]
       (\sqrt{\sum_i y_i\l_{is}} - \sqrt{\sum_i y_i\l_{im}})
      = \displaystyle {y_i\over 2}\Big({1\over\sqrt{\sum_i y_i\l_{is}}}-{1\over\sqrt{\sum_i y_i\l_{im}}}\Big)
      \leq 0.
\eqa
That is, if we decrease
$\l_{is}\leadsto\l'_{is}:=\l_{is}-\delta_i$ and
$\l_{im}\leadsto\l'_{im}:=\l_{im}-\delta_i$ by the same amount
$\delta_i$, then (\ref{eqsqrtl}) increases. The maximal possible
$\delta_i:=\min\{\l_{is},\l_{im}\}$ makes $\l'_{is}$ or $\l'_{im}$
zero, hence $0\leq\l'_{is}+\l'_{im}\leq 1$. Similarly
\beqn\textstyle
  0 \;\leq\; \sqrt{\sum_i z_i\l_{im}} - \sqrt{\sum_i z_i\l_{is}}
    \;\leq\; \sqrt{\sum_i z_i\l'_{im}} - \sqrt{\sum_i z_i\l'_{is}}
\eeqn
This implies
\bqan
  0 &\leq& \textstyle \sqrt{l_t^{\smash{\Lambda_\xi}}} - \sqrt{l_t^{\smash{\Lambda_\mu}}}
  \;\equiv\; \sqrt{\sum_i y_i\l_{is}} - \sqrt{\sum_i y_i\l_{im}}
\\
  &\leq& \textstyle \sqrt{\sum_i y_i\l'_{is}} - \sqrt{\sum_i y_i\l'_{im}}
      \;+\;  \sqrt{\sum_i z_i\l'_{im}} - \sqrt{\sum_i z_i\l'_{is}}
\\
  &\leq& \textstyle \sqrt{\sum_i\l'_{is}(\sqrt{y_i}\!-\!\sqrt{z_i})^2}
              +     \sqrt{\sum_i\l'_{im}(\sqrt{y_i}\!-\!\sqrt{z_i})^2}
\\
  &\leq& \textstyle \sqrt{2\sum_i(\l'_{is}\!+\!\l'_{im})(\sqrt{y_i}\!-\!\sqrt{z_i})^2}
   \;\leq\; \sqrt{2\sum_i(\sqrt{y_i}\!-\!\sqrt{z_i})^2}
   \;\equiv\; \sqrt{2h_t}
\eqan
In the third inequality I used the Hellinger bound
\req{eqHellBnd} twice, and in the fourth inequality I used
$\sqrt{a}+\sqrt{b}\leq\sqrt{2(a+b)}$.
Without the reduction
$\l\leadsto\l'$ the bound would have been a factor of $\sqrt{2}$
worse.
Taking the square, expectation, and sum over $t$ proves the last
inequality in \req{lbnd}.
The first inequality in \req{lbnd} is again an
instantiation of \req{eqHellBnd} with $i\leadsto (t,\o_{<t})$
and $v_i\leadsto\mu(\o_{<t})$, i.e.\ $\sum_i v_i...\leadsto
\sum_t\E[...]$ and $a_i\leadsto l_t^{\smash{\Lambda_\xi}}$
and $b_i\leadsto l_t^{\smash{\Lambda_\mu}}$. \qed

\paradot{Proof of equation (\ref{Kprop}$f$)}
Function $P(y):=\sum_{x:f(x)=y} 2^{-K(x)}$ is lower
semi-computable, since $K(x)$ is upper semi-computable, all $x\in\X^*$
can be enumerated, and $f(x)=y$ is decidable. Further,
$\sum_y P(y)=\sum_x 2^{-K(x)}\leq 1$, hence
MDL bound (\ref{Kprop}$e$) implies $K(y)\leqa -\lb P(y)+K(P)$.
Let $g(y)=\min\{x:f(y)=x\}$ be the lexicographically first inverse of $f$.
With $K(P)\leqa K(f)=O(1)$, also function $g$ has complexity $O(1)$.
Hence
\beqn
  2^{-K(y)} \;\geqm\; P(y)
  \;\equiv\; \sum_{x:f(x)=y} 2^{-K(x)} \;\geq\; 2^{-K(g(y))}
  \;\geqm\; 2^{-K(y)}
\eeqn
where I dropped all but the contribution from $g(y)$ in the sum,
and used (\ref{Kprop}$d$) for $g$. \qed

\paradot{Proof of bound \req{iMbnd} $M(\bar x_n|x_{<n})\geqm 2^{-K(n)}$}
For $x=x_{<n}\in\X^{n-1}$ and $a=\bar x_n\in\X$ we have
\ifjournal
\beqn
  M(a|x) \stackrel{(a)}=  {M(xa)\over M(x)}
  \stackrel{(b)}= { \sum_{p:U(p)=xa*} 2^{-\l(p)} \over \sum_{p:U(p)=x*} 2^{-\l(p)} }
  \stackrel{(c)}\geq { \sum_{p:U(\tilde p)=xa} 2^{-\l(\tilde p)} \over \sum_{p:U(p)=x*} 2^{-\l(p)} }
  \stackrel{(d)}= 2^{-\l(qn^*)}
  \stackrel{(e)}{\stackrel\times=} 2^{-K(n)}
\eeqn
\else
\beqn
  M(a|x) \;\stackrel{(a)}=\;  {M(xa)\over M(x)}
  \;\stackrel{(b)}=\; { \sum_{p:U(p)=xa*} 2^{-\l(p)} \over \sum_{p:U(p)=x*} 2^{-\l(p)} }
  \;\stackrel{(c)}\geq\; { \sum_{p:U(\tilde p)=xa} 2^{-\l(\tilde p)} \over \sum_{p:U(p)=x*} 2^{-\l(p)} }
  \;\stackrel{(d)}=\; 2^{-\l(qn^*)}
  \;\stackrel{(e)}{\stackrel\times=}\; 2^{-K(n)}
\eeqn
\fi
In $(a)$ and $(b)$ I simply inserted the definition \req{Mdef}
of $M$.
I now $(c)$ restrict the sum over all $p:U(p)=xa*$ in the
numerator to programs $\tilde p$ of the following form: $\tilde
p=q n^*p$, where $U(p)=x*$, $n^*$ is the shortest code of $n$, and
$q$ simulates $p$ until $n-1$ symbols are printed, then prints $a$,
and thereafter halts, i.e.\ $U(\tilde p)=xa$.
The numerator now sums over exactly the same programs $p$ as the
denominator. Since $2^{-\l(\tilde
p)}=2^{-\l(q n^*)}2^{-\l(p)}$, and $2^{-\l(q n^*)}$ is a constant
independent of $p$, numerator and denominator cancel and $(d)$ follows.
$(e)$ follows from the definition of $n^*$ and from $\l(q)=O(1)$.
\qed

\paradot{Proof of bound \req{iMbnd} $M(\bar x_n|x_{<n})\leqm 2^{2\Km(x_{1:n})-K(n)}$}
Assume $x_{1:\infty}$ is a computable sequence, $\X$ is binary, and
$\bar x_n\neq x_n$, and define $P(n):=M(x_{<n}\bar x_n)$. Given
$x_{1:\infty}$, $P$ can be semi-computed from below, hence
$K(P)\leqa\Km(x_{1:\infty})$. Also $\sum_n P(n)\leq 1$, since $\{
x_{<n}\bar x_n : n\in\SetN\}$ forms a prefix-free set. Hence
$K(n)\leqa-\lb P(n)+K(P)$ by (\ref{Kprop}$e$), which implies
$M(x_{<n}\bar x_n)\leqm 2^{\Km(x_{1:\infty})-K(n)}$. Since
$M(x_{<n})\geq 2^{-\Km(x_{<n})}\geq 2^{-\Km(x_{1:\infty})}$, we get
$M(\bar x_n|x_{<n}) \leqm 2^{2\Km(x_{1:\infty})-K(n)}$, which nearly
is \req{iMbnd}. Since the l.h.s.\ is independent of
$x_{n+1:\infty}$, a bound independent of it should be (and is)
possible, as we will now show.

Consider sequence $x_{1:n}$ and shortest program $p$ printing
$x_{1:n}*$. Let $U_t$ be $U$ stopped after $t$ time steps and
define corresponding $M_t$. Then $U_t(p)=x_{1:n_t}$ (for some
$x_{n+1:n_t}$ if $n_t>n$). I define $P_t(n'):=\sum_{a\neq
x_{n'}}M_t(x_{<n'}a)$ for $n'\leq n_t$ and 0 for $n'>n_t$. With $n_t$
also $P_t$ is computable and increasing, hence
$P(n'):=\lim_{t\to\infty} P_t(n')=\sup_t P_t(n')$ is lower
semi-computable. Clearly $P(n')=\sum_{a\neq x_{n'}}M(x_{<n'}a)$ for
$n'\leq n_\infty$ and $P(n')=0$ for $n'>n_\infty$
($n'_\infty=\lim_t n_t\leq\infty$). Hence $\sum_{n'}P(n')\leq 1$, since
$\{x_{<n'}a\,:\,a\neq x_{n'},\,n'\leq n_\infty\}$ is a prefix free
set, which implies $K(n)\leqa-\lb P(n)+K(P)$ by (\ref{Kprop}$e$).
Since $n\leq n_\infty$ and $K(P)\leqa \l(p)=\Km(x_{1:n})$, we get
$\sum_{a\neq x_n}M(x_{<n}a)\leqm 2^{\Km(x_{1:n})-K(n)}$. Using
$M(x_{<n})\geq 2^{-\Km(x_{<n})}\geq 2^{-\Km(x_{1:n})}$, we get the
desired bound $M(\bar x_n|x_{<n}) \leq \sum_{a\neq x_n}
M(a|x_{<n}) \leqm 2^{2\Km(x_{1:n})-K(n)}$. \qed

\end{document}